\documentclass[11pt]{article}
\usepackage{amsmath,amssymb,amsthm,color}

\usepackage{comment}
\usepackage{color}
\usepackage{amssymb}
\usepackage{amsmath}
\usepackage[english]{babel}
\usepackage{fancyhdr}
\usepackage{amsmath}

\def\ci{\begin{color}{red}\,}
\def\cf{\end{color}\,}

\newtheorem{thm}{Theorem}
\newtheorem{corollary}[thm]{Corollary}

\newtheorem{claim}[thm]{Claim}

\newtheorem{theorem}[thm]{Theorem}
\newtheorem{remark}[thm]{Remark}

\begin{document}
\begin{center}
{\bf Nonlinear $*$-Jordan-type derivations on alternative $*$-algebras}
\\
\vspace{.2in}

{\bf Aline Jaqueline de Oliveira Andrade}
\\
Federal University of ABC\\
dos Estados Avenue, 5001\\
09210-580, Santo Andr\'{e}, Brazil.
\\
aline.jaqueline@ufabc.edu.br
\\
\vspace{.2in}

{\bf Gabriela C. Moraes}
\\
Federal University of ABC\\
dos Estados Avenue, 5001\\
09210-580, Santo Andr\'{e}, Brazil.
\\
gcotrim@gmail.com
\\
\vspace{.2in}

{\bf Ruth Nascimento Ferreira}
\\
Federal University of Technology,\\
Professora Laura Pacheco Bastos Avenue, 800\\
85053-510, Guarapuava, Brazil.
\\
ruthnascimento@utfpr.edu.br
\\
\vspace{.2in}

{\bf Bruno Leonardo Macedo Ferreira}
\\
Federal University of Technology,\\
Professora Laura Pacheco Bastos Avenue, 800,\\
85053-510, Guarapuava, Brazil.
\\
brunoferreira@utfpr.edu.br
\\
brunolmfalg@gmail.com
\\

\end{center}

\ 


\

AMS: 17D05, 47B47.

\begin{abstract}
Let $A$ be an unital alternative $*$-algebra. Assume that $A$ contains
a nontrivial symmetric idempotent element $e$ which satisfies $xA \cdot e = 0$ implies $x = 0$ and
$xA \cdot (1_A - e) = 0$ implies $x = 0$. In this paper, it is shown that $\Phi$ is a
nonlinear $*$-Jordan-type derivation on A if and only if $\Phi$ is an additive
$*$-derivation. As application, we get a result on alternative $W^{*}$-algebras. 
\end{abstract}

\vspace{0.5 in}

\section*{Introduction}

An algebra $A$ not necessarily associative or commutative is called an alternative
algebra if it satisfies the identities $a^2b = a(ab)$ and $ba^2 = (ba)a$, for all elements
$a, b \in A$. One easily sees that any associative algebra is an alternative algebra. An
alternative algebra $A$ is called prime if for any elements $a, b \in A$ satisfying the
condition $aAb = 0$, then either $a = 0$ or $b = 0$.
Let $A$ and $B$ be alternative algebras. We say that a mapping $\Phi : A \rightarrow B$
preserves product if $\Phi(ab) = \Phi(a)\Phi(b)$, for all elements $a, b \in A$, and that it preserves Jordan product (resp., preserves Lie product) if $\Phi(ab + ba) = \Phi(a)\Phi(b) + \Phi(b)\Phi(a)$ (resp., $\Phi(ab - ba) = \Phi(a)\Phi(b) - \Phi(b)\Phi(a))$, for all elements $a, b \in A$.
We say that a mapping $\Phi : A \rightarrow B$ is additive if $\Phi(a + b) = \Phi(a) + \Phi(b)$,
for all elements $a, b \in A$ and that it is an isomorphism if $\Phi$ is an additive bijection that preserves products.
By involution, we mean a mapping $* : A \rightarrow A$ such that $(x + y)^{*} = x^{*} + y^{*}; (x^{*})^{*} = x$ and $(xy)^{*} = y^{*}x^{*}$ for all $x, y \in A$. An
element $s \in A$ satisfying $s^{*} = s$ is called a symmetric element of $A$.
Let $A$ and $B$ be alternative $*$-algebras. We say that a mapping $\Phi : A \rightarrow B$ preserves involution if $\Phi(a^{*}) = \Phi(a)^{*}$, for all elements $a \in A$, and that $\Phi$ is a $*$-isomorphism if $\Phi$ is an isomorphism that preserves involution.

Let $A$ be a alternative $*$-algebra over the ground field $\mathbb{C}$. For any $a,b \in A$ denote a new
product of $a$ and $b$ by $a\bullet b = ab+ba^{*}$, and this new product $\bullet$ is usually known
as Jordan $\bullet$-product. Such kind of product plays a more and more important role
in some research topics, and its study has attracted many authors' attention.
Particular attention has been paid to understanding maps which preserve the
product $AB + BA^{*}$ on $*$-algebra (for example, see \cite{Dai, Huo, LiLuFa, LiLu, LiLu2, ZhaLi}). The product
is extensively studied because it naturally arise in the problem of representing quadratic functionals with sesquilinear functionals (for example, see \cite{Sem1, Sem3, Sem3}) and
in the problem of characterizing ideals (for example, see \cite{BreFo, Molnar}).  Furthermore, there has been a great interest in the study of additivity or characterization 
of mappings that preserve some kind of product on  non-associative rings  and algebras for example in the paper`s  \cite{Fer,fgk20,fgf20,f19,mff16}.
For these kinds of mappings defined above, in a recent paper Ferreira and Costa \cite{FerreiraCosta} has studied the characterization about $*$-Jordan-type maps on $C^{*}$-algebras, they have proved under mild conditions that every multiplicative $*$-Jordan-type maps are $*$-isomorphism. That work motivated us to study the same question for the case $*$-Jordan-type derivations on a class of nonassociative algebras, namely $*$-alternative algebras. Let us define $*$-Jordan-type derivations.

Recall that an additive map $\Phi : A \rightarrow A$ is said to be an additive derivation
if $\Phi(ab) = \Phi(a)b + a\Phi(b)$ for all $a,b \in A$ Furthermore, $\Phi$ is said to be an
additive $*$-derivation if it is an additive derivation and satisfies $\Phi(a^*) = \Phi(a)^*$
for all $a \in A$. We say that $\Xi_{y,z}: A \rightarrow A$ is an additive $*$-inner derivation if it is an additive inner derivation, that is,
$$\Xi_{y,z}(a) = ([L_{y},L_{z}] + [L_y, R_z] + [R_y,R_z])(a),$$
is holds true for all $y, z \in A$ and satisfies $\Xi(a^{*}) = \Xi(a)^{*}$. A not necessarily linear map $\Phi : A \rightarrow A$ is said to be a nonlinear
$*$-Jordan derivation if
$$\Phi(a \bullet b) = \Phi(a) \bullet b + a \bullet \Phi(b)$$
for all $a,b \in A$. Similarly, a map $\Phi : A \rightarrow A$ is said to be a nonlinear $*$-Jordan
triple derivation if
$$\Phi(a \bullet b \bullet c) = \Phi(a) \bullet b \bullet c + a \bullet \Phi(b) \bullet c + a \bullet b \bullet \Phi(c)$$
for all $a,b,c \in A$, where $a \bullet b \bullet c = (a \bullet b) \bullet c$ (we should be aware that $\bullet$
is not necessarily associative).
Given the consideration of nonlinear $*$-Jordan derivations and nonlinear $*$-
Jordan triple derivations, as well as introduced the definition of $*$-Jordan-type maps in \cite{FerreiraCosta} we introduce the following. Suppose that $n \geq 2$ is a fixed positive integer. Accordingly, a nonlinear $*$-Jordan $n$-derivation is a map $\Phi : A \rightarrow A$ satisfying the condition
\begin{eqnarray}\label{ouro}
\Phi(a_1 \bullet a_2 \bullet \cdots \bullet a_n) =
\sum_{k=1}^n a_1 \bullet \cdots \bullet a_{k-1} \bullet \Phi(a_k) \bullet a_{k+1} \bullet \cdots \bullet a_n
\end{eqnarray}
for all $a_1, a_2, \ldots a_n \in A$, where $a_1 \bullet a_2 \bullet \cdots \bullet a_n = (\cdots ((a_1 \bullet a_2)\bullet a_3) \cdots \bullet a_n)$. 
By the definition, it is clear that every $*$-Jordan derivation is a $*$-Jordan $2$-derivation and every $*$-Jordan triple derivation is a $*$-Jordan $3$-derivation.
It is obvious that every nonlinear $*$-Jordan derivation on any $*$-algebra is a $*$-Jordan
$n$-derivation. But we do not know whether the converse is true. $*$-Jordan $2$-derivation, $*$-Jordan 3-derivation and $*$-Jordan n-derivation are collectively referred to as $*$-Jordan-type derivations. $*$-Jordan-type derivations in different
backgrounds are extensively studied by several authors, see \cite{Darvish, LiLuFang2, TaghaRo, Zhang, ZhaoLi}.

It is well known that $A$ has a Peirce decomposition
$A=A_{11}\oplus A_{12}\oplus A_{21}\oplus A_{22},$ where
$A_{ij}=e_{i}Ae_{j}$ $(i,j=1,2)$ \cite{He}, satisfying the following multiplicative relations:
\begin{enumerate}\label{asquatro}
\item [\it (i)] $A_{ij}A_{jl}\subseteq A_{il}\
(i,j,l=1,2);$
\item [\it (ii)] $A_{ij}A_{ij}\subseteq A_{ji}\
(i,j=1,2);$
\item [\it (iii)] $A_{ij} A_{kl}=0,$ if $j\neq k$ and
$(i,j)\neq (k,l),\ (i,j,k,l=1,2);$
\item [\it (iv)] $x_{ij}^{2}=0,$ for all $x_{ij}\in A_{ij}\ (i,j=1,2;~i\neq j).$
\end{enumerate}

The notion of Peirce decomposition for alternative algebras is similar to that one for associative algebras. However, this similarity is restricted to its written form, not including its theoretical structure since Peirce decomposition for alternative algebras is a generalization of that classical one for associative algebras.

\section{Main theorem}

We shall prove as follows the main result of this paper.

\begin{theorem}\label{mainthm} 
Let $A$ be a unital alternative $*$-algebra with the unit $1_A$. Assume that $A$
contains a nontrivial symmetric idempotent element $e$ which satisfies
\begin{enumerate}
\item[$(\spadesuit)$] $xA \cdot e = 0$ implies $x = 0$;
\item [$(\clubsuit)$] $xA \cdot (1_A - e) = 0$ implies $x = 0$.
\end{enumerate}
If a map $\Phi : A \rightarrow A$ satisfies
$$\Phi(a_1 \bullet a_2 \bullet \cdots \bullet a_n) =
\sum_{k=1}^n a_1 \bullet \cdots \bullet a_{k-1} \bullet \Phi(a_k) \bullet a_{k+1} \bullet \cdots \bullet a_n$$
for all $a_{n-1}, a_n \in A$ and $a_i = 1_A$ for all $i \in \left\{1, 2, \ldots, n-2\right\}$, then $\Phi$ is an additive $*$-derivation.  
\end{theorem}

\begin{corollary}
Let $A$ be a unital alternative $*$-algebra with the unit $1_A$. Assume that $A$
contains a nontrivial symmetric idempotent element $e$ which satisfies $(\spadesuit)$ and $(\clubsuit)$.
Then $\Phi$ is a nonlinear $*$-Jordan-type derivation on $A$ if and only if $\Phi$ is an additive
$*$-derivation.
  
\end{corollary}

It is obvious that prime alternative $*$-algebras satisfy the assumptions ``$x A \cdot e = 0 \Rightarrow x = 0$ and $x A \cdot (1_A-e) = 0 \Rightarrow x = 0$". Thus, we have the following result

\begin{corollary}
Let $A$ be a prime unital alternative $*$-algebra with the unit $1_A$. Assume that $A$
contains a nontrivial symmetric idempotent element $e$. If a map $\Phi : A \rightarrow A$ satisfies
$$\Phi(a_1 \bullet a_2 \bullet \cdots \bullet a_n) =
\sum_{k=1}^n a_1 \bullet \cdots \bullet a_{k-1} \bullet \Phi(a_k) \bullet a_{k+1} \bullet \cdots \bullet a_n$$
for all $a_{n-1}, a_n \in A$ and $a_i = 1_A$ for all $i \in \left\{1, 2, \ldots, n-2\right\}$, then $\Phi$ is an additive $*$-derivation.
\end{corollary}

 \begin{remark}
Some examples of non-prime alternative algebras satisfying  the assumptions ``$x A \cdot e = 0 \Rightarrow x = 0$ and $x A \cdot (1_A - e) = 0 \Rightarrow x = 0$"
were given in \cite{posd}. 
\end{remark}

\vspace{.1in}
\section{The proof of main result}

Consider the Peirce decomposition of $A$ to respect nontrivial symmetric idempotent element $e$. In all that follows, when we write $a_{ij}$, it indicates that $a_{ij} \in A_{ij}$ and $e_1 = e$, $e_2 = 1_A - e_1$. We will complete the proof by proving several claims.

Borrowing the notation in \cite{fgf20} let us point out a non-associative monomial $m$ of degree $n$ by the following expression form
\[
m(x_1, x_2, \cdots, x_n, x_{n-1}) := ( ( \cdots (x_1 \bullet x_2) \cdots)\bullet x_n).
\]
When the first $n-2$ variables in the non-associative monomial $m$ assume equal values, we
will denote by
\begin{gather*}
    m(\underbrace{a, \cdots, a}_{n-2}, x_{n-1}, x_n) := \xi_{a}(x_{n-1}, x_n), \ \textup{ for any} \ a \in A.
\end{gather*}

\begin{claim}\label{claim1}
$\Phi(0) = 0.$
\end{claim}

Indeed, 
\[
\Phi(0) = \Phi(\xi_{1_A}(0, 0)) = 0.
\]

\begin{claim}\label{claim2}
For every $a_{12} \in A_{12}$, $b_{21} \in A_{21}$ we have
$$\Phi(a_{12} + b_{21}) = \Phi(a_{12}) + \Phi(b_{21}).$$

\end{claim}
Consider $t = \Phi(a_{12} + b_{21}) - \Phi(a_{12}) - \Phi(b_{21})$.
Since
\[
\xi_{1_A}((e_1-e_2), a_{12}) = \xi_{1_A}((e_1-e_2), b_{21}) = 0.
\]

It follows from Claim \ref{claim1} that
\begin{align*}
\sum_{\Phi(1_A)} & \xi_{\Phi(1_A)}(e_1-e_2, a_{12}+b_{21}) + \xi_{1_A}( \Phi(e_1 - e_2), (a_{12} + b_{21})) \\
& \ + \xi_{1_A}((e_1 - e_2), \Phi(a_{12} + b_{21})) \\
=& \  \Phi( \xi_{1_A}(e_1-e_2, a_{12}+b_{21})) \\
=& \  \Phi( \xi_{1_A}(e_1-e_2, a_{12})) + \Phi(\xi_{1_A}(e_1-e_2, b_{21})) \\
=& \ \sum_{\Phi(1_A)} \xi_{\Phi(1_A)}(e_1-e_2, a_{12}+b_{21}) + \xi_{1_A}( \Phi(e_1 - e_2), (a_{12} + b_{21}))\\ & \ + \xi_{1_A}(e_1-e_2, \Phi(a_{12}) + \Phi(b_{21})). \end{align*}

From this, we get that $\xi_{1_A}(e_1-e_2, t)=0$ 
Since $\xi_{1_A}(a_{12},e_1)=0,$ we have 
\[
\xi_{1_A}(a_{12}+b_{21}, e_1) = \xi_{1_A}(b_{21}, e_1).
\]
It follows that
\begin{align*}
\sum_{\Phi(1_A)} & \xi_{\Phi(1_A)}(a_{12}+b_{21}, e_1) + \xi_{1_A}\left( \Phi(a_{12}+b_{21}, e_1)\right) + \xi_{1_A}(a_{12}+b_{21},\Phi(e_1)) \\
=& \ \Phi(\xi_{1_A}(a_{12}+b_{21},e_1)) \\
=& \ \Phi(\xi_{1_A}(a_{12}, e_1)) + \Phi(\xi_{1_A}(b_{21},e_1))\\
=& \ \sum_{\Phi(1_A)} \xi_{\Phi(1_A)}(a_{12}+b_{21}, e_1) + \xi_{1_A} (\Phi(a_{12})+\Phi(b_{21}), e_1) + \xi_{1_A} (a_{12}+b_{21},e_1).
\end{align*}
Hence $\xi_{1_A}(t,e_1)=0$ from which we get that $t_{21} = 0$. Similarly, $t_{12} = 0$,
proving the claim. 

\begin{claim}\label{claim3}
For every $a_{11} \in A_{11}$, $b_{12} \in A_{12}$, $c_{21} \in A_{21}$, $d_{22} \in A_{22}$ we have
$$\Phi(a_{11} + b_{12} + c_{21}) = \Phi(a_{11}) + \Phi(b_{12}) + \Phi(c_{21})$$
and
$$\Phi(b_{12} + c_{21} + d_{22}) = \Phi(b_{12}) + \Phi(c_{21}) + \Phi(d_{22}).$$
\end{claim}

Let $t = \Phi(a_{11} + b_{12} + c_{21}) - \Phi(a_{11}) - \Phi(b_{12}) - \Phi(c_{21})$. Since $\xi_{1_A}(e_2,a_{11})=0$ it follows from Claim \ref{claim2} that 
\begin{align*}
\sum_{\Phi(1_A)} & \xi_{\Phi(1_A)} (e_2, a_{11}+b_{12}+c_{21}) + \xi_{1_A}(\Phi(e_2),a_{11}+b_{12}+c_{21}) + \xi_{1_A} (e_2, \Phi(a_{11}+b_{12}+c_{21})) \\
=& \ \Phi(\xi_{1_A}(e_2, a_{11}+b_{12}+c_{21})) \\
=& \ \Phi(\xi_{1_A}(e_2, a_{11})) + \Phi(\xi_{1_A}(e_2,b_{12}+c_{21})) \\
=& \ \Phi(\xi_{1_A}(e_2, a_{11})) + \Phi(\xi_{1_A}(e_2, b_{12})) + \Phi(\xi_{1_A}(e_2, c_{21}))\\
=& \ \sum_{\Phi(1_A)} \xi_{\Phi(1_A)}(e_2,a_{11}+b_{12}+c_{21}) + \xi_{1_A} (\Phi(e_2), a_{11}+b_{12}+c_{21}) \\ 
& \ + \xi_{1_A}(e_2, \Phi(a_{11}) + \Phi(b_{12})+\Phi(c_{21}))
\end{align*}
Hence $\xi_{1_A}(e_2,t)=0$ from which we get that $t_{12} = t_{21} = t_{22} = 0$. 
Since
\[
\xi_{1_A}(e_1-e_2, a_{12}) = \xi_{1_A}(e_1-e_2, a_{21})=0
\]
one has
\begin{align*}
\sum_{\Phi(1_A)} & \xi_{\Phi(1_A)} (e_1-e_2,a_{11}+b_{12}+c_{21}) + \xi_{1_A}(\Phi(e_1-e_2), a_{11}+b_{12}+c_{21}) \\
& \ + \xi_{1_A}(e_1-e_2, \Phi(a_{11}+b_{12}+c_{21}))\\
=& \ \Phi(\xi_{1_A}(e_1-e_2, a_{11}+b_{12}+c_{21})) \\
=& \ \Phi(\xi_{1_A}(e_1-e_2, a_{11})) + \Phi(\xi_{1_A}(e_1-e_2, b_{12})) + \Phi(\xi_{1_A}(e_1-e_2, c_{21}))\\
=& \ \sum_{\Phi(1_A)} \xi_{\Phi(1_A)} (e_1-e_2, a_{11}+b_{12}+c_{21}) + \xi_{1_A}(\Phi(e_1-e_2),a_{11}+b_{12}+c_{21})\\
& \ + \xi_{1_A}(e_1-e_2, \Phi(a_{11})+ \Phi(b_{12}) + \Phi(c_{21}))
\end{align*}
from which we get that $\xi_{1_A}(e_1-e_2,t)=0.$ Thus, $t_{11} = 0$, and then $t = 0$. 
Similarly, we can prove that $\Phi(b_{12} + c_{21} + d_{22}) = \Phi(b_{12}) + \Phi(c_{21}) + \Phi(d_{22}).$

\begin{claim}\label{claim4}
For every $a_{11} \in A_{11}$, $b_{12} \in A_{12}$, $c_{21} \in A_{21}$, $d_{22} \in A_{22}$ we have
$$\Phi(a_{11} + b_{12} + c_{21} + d_{22}) = \Phi(a_{11}) + \Phi(b_{12}) + \Phi(c_{21}) + \Phi(d_{22}).$$

\end{claim}

Let $t = \Phi(a_{11} + b_{12} + c_{21} + d_{22}) - \Phi(a_{11}) - \Phi(b_{12}) - \Phi(c_{21}) - \Phi(d_{22}).$
Since $\xi_{1_A}(e_2,a_{11})=0$ it follows from Claim \ref{claim3} that 
\begin{align*}
\sum_{\Phi(1_A)} & \xi_{\Phi(1_A)} (e_2, a_{11} + b_{12}+c_{21} + d_{22}) + \xi_{1_A}(\Phi(e_2, a_{11} + b_{12}+c_{21} + d_{22}))\\
& \ + \xi_{1_A}(e_2, \Phi(a_{11} + b_{12}+c_{21} + d_{22})) \\
=& \ \Phi(\xi_{1_A}(e_2, a_{11} + b_{12}+c_{21} + d_{22})) \\
=& \ \Phi(\xi_{1_A}(e_2, a_{11})) + \Phi(\xi_{1_A}(e_2,b_{12}+c_{21} + d_{22}))\\
=& \ \Phi(\xi_{1_A}(e_2, a_{11})) + \Phi(\xi_{1_A}(e_2, b_{12})) + \Phi(\xi_{1_A}(e_2, c_{21})) + \Phi(\xi_{1_A}(e_2, d_{22}))\\
=& \ \sum_{\Phi(1_A)} \xi_{\Phi(1_A)} (e_2, a_{11} + b_{12}+c_{21} + d_{22}) + \xi_{1_A} (\Phi(e_2), a_{11} + b_{12}+c_{21} + d_{22}) \\
& \ + \xi_{1_A} (e_2, \Phi(a_{11})+\Phi(b_{12})+\Phi(c_{21})+\Phi(d_{22}))
\end{align*}

Hence $\xi_{1_A}(e_2,t)=0$ from which we get that $t_{12} = t_{21} = t_{22} = 0$. 
Similarly, we can prove $t_{11} = 0$ proving the claim.

In the next Claim we have a slight difference from the case of associative algebras. In the associative case we know that $a_{ij}b_{ij} = 0$ for $i \neq j$ but in general in alternative algebras we do not have this property.

\begin{claim}\label{claim5}
For every $a_{12}, b_{12} \in A_{12}$ and $c_{21}, d_{21} \in A_{21}$ we have  
$$\Phi(a_{12}b_{12}+a^{*}_{
12}) = \Phi(a_{12}b_{12})+\Phi(a^{*}_{
12})$$
and
$$\Phi(c_{21}d_{21}+c^{*}_{
21}) = \Phi(c_{21}d_{21})+\Phi(c^{*}_{
21}).$$
\end{claim}
Since 
\[\xi_{1_A}\left( a_{12}, \dfrac{e_2+b_{12}}{2n-2} \right) = a_{12}+a_{12}b_{12} + a_{12}^*+b_{12}a_{12}^*
\]
we get from Claim \ref{claim4} that
\begin{align*}
\Phi(a_{12}) & + \Phi(a_{12}b_{12} + a^{*}_{12}) + \Phi(b_{12}a^{*}_{12}) \\
=& \ \Phi(a_{12} + a_{12}b_{12} + a^{*}_{12} + b_{12}a^{*}_{12}) \\
=& \ \Phi\left( \xi_{1_A}\left( a_{12}, \dfrac{e_2+b_{12}}{2n-2} \right) \right) \\
=& \ \sum_{\Phi(1_A)} \xi_{\Phi(1_A)} \left( a_{12}, \dfrac{e_2+b_{12}}{2n-2} \right) + \xi_{1_A} \left( \Phi(a_{12}), \dfrac{e_2+b_{12}}{2n-2} \right) + \xi_{1_A}\left( a_{12}, \Phi \left( \dfrac{e_2+b_{12}}{2n-2} \right) \right) \\
=& \ \sum_{\Phi(1_A)} \xi_{\Phi(1_A)} \left( a_{12}, \dfrac{e_2+b_{12}}{2n-2} \right) + \xi_{1_A} \left( \Phi(a_{12}), \dfrac{e_2+b_{12}}{2n-2} \right) \\ 
& \ + \xi_{1_A} \left( a_{12}, \Phi \left( \dfrac{e_2}{2n-2} \right) + \Phi \left( \dfrac{b_{12}}{2n-2} \right) \right) \\
=& \ \Phi \left( \xi_{1_A} \left( a_{12}, \dfrac{e_2}{2n-2} \right) \right) + \Phi \left( \xi_{1_A} \left( a_{12}, \dfrac{b_{12}}{2n-2} \right) \right) \\
=& \ \Phi(a_{12}) + \Phi(a_{12}^*) + \Phi(a_{12}b_{12}) + \Phi(b_{12}a_{12}^*),
\end{align*}
which implies $\Phi(a_{12}b_{12} + a^{*}_{12}) = \Phi(a_{12}b_{12}) + \Phi(a^{*}_{12}).$
Similarly, we prove the other case using the identity
\[
\xi_{1_A}\left( c_{21}, \dfrac{e_1+d_{21}}{2n-2} \right) = c_{21}+ c_{21}d_{21}+c_{21}^* + d_{21}c_{21}^*.
\]
\begin{claim}\label{claim6}
For all $a_{ij}, b_{ij} \in A_{ij}$, $1 \leq i \neq j \leq 2$, we have
$$\Phi(a_{ij} + b_{ij}) = \Phi(a_{ij}) + \Phi(b_{ij}).$$
\end{claim}

Since
\[
\xi_{1_A} \left( \dfrac{e_i+a_{ij}}{2n-2}, e_j+b_{ij} \right) = a_{ij} + b_{ij} + a_{ij}^* + a_{ij}b_{ji}+b_{ij}a_{ij}^*
\]
we get from Claim \ref{claim5} that
\begin{align*}
\Phi & (a_{ij}+b_{ij}) \Phi(a_{ij}^* + a_{ij}b_{ij}) + \Phi(b_{ij}a_{ij}^*) \\
=& \ \Phi \left( \xi_{1_A} \left( \dfrac{e_i + a_{ij}}{2n-2}, e_j+b_{ij} \right) \right)\\
=& \ \sum_{\Phi(1_A)} \xi_{\Phi(1_A)} \left( \dfrac{e_i + a_{ij}}{2n-2}, e_j+b_{ij} \right) + \xi_{1_A} \left( \Phi \left( \dfrac{e_i + a_{ij}}{2n-2}\right), e_j+b_{ij}\right)\\
& \ + \xi_{1_A} \left( \dfrac{e_i + a_{ij}}{2n-2}, \Phi(e_j+b_{ij}) \right)\\
=& \ \sum_{\Phi(1_A)} \xi_{\Phi(1_A)}\left( \dfrac{e_i}{2n-2} + \dfrac{a_{ij}}{2n-2}, e_j+b_{ij} \right) + \xi_{1_A} \left( \Phi \left( \dfrac{e_i}{2n-2} \right) + \Phi\left( \dfrac{a_{ij}}{2n-2} \right), e_i + b_{ij}\right)\\
& \ + \xi_{1_A} \left(  \dfrac{e_i}{2n-2} + \dfrac{a_{ij}}{2n-2}, \Phi(e_j)+\Phi(b_{ij}) \right)\\
=& \ \Phi \left( \xi_{1_A} \left(\dfrac{e_i}{2n-2},e_j \right)\right)+\Phi \left( \xi_{1_A} \left(\dfrac{a_{ij}}{2n-2},b_{ij} \right)\right) + \Phi\left( \xi_{1_A} \left( \dfrac{a_{ij}}{2n-2}, e_j \right) \right) \\
& \ + \Phi \left( \xi_{1_A} \left( \dfrac{e_i}{2n-2}, b_{ij} \right) \right) \\
=& \ \Phi(b_{ij}) + \Phi(a_{ij} + a_{ij}^*) + \Phi(a_{ij}b_{ij}+b_{ij}a_{ij}^*)\\ 
=& \Phi(b_{ij}) + \Phi(a_{ij}) + \Phi(a_{ij}^*) + \Phi(a_{ij}b_{ij}) + \Phi(b_{ij}a_{ij}^*), 
\end{align*}
which implies $\Phi(a_{ij} + b_{ij}) = \Phi(b_{ij}) + \Phi(a_{ij})$.

\begin{claim}\label{claim7}
For every $a_{ii}, b_{ii} \in A_{ii}$, $1 \leq i \leq 2$ we have
$$\Phi(a_{ii} + b_{ii}) = \Phi(a_{ii}) + \Phi(b_{ii}).$$
\end{claim}
Let $t = \Phi(a_{ii} + b_{ii}) - \Phi(a_{ii}) - \Phi(b_{ii})$. 
Since 
\[
\xi_{1_A}(e_j, a_{ii}) = \xi_{1_A}(e_j,b_{ii})=0,
\]
it follows that
\begin{align*}
\sum_{\Phi(1_A)} & \xi_{\Phi(1_A)} (e_j, a_{ii}+b_{ii}) + \xi_{1_A} (\Phi(e_j), a_{ii}+b_{ii}) + \xi_{1_A}(e_j, \Phi(a_{ii}+b_{ii}))\\
=& \ \Phi(\xi_{1_A}(e_j, a_{ii}+b_{ii}))\\
=& \ \Phi(\xi_{1_A}(e_j, a_{ii})) + \Phi(\xi_{1_A}(e_j, b_{ii}))\\
=& \ \sum_{\Phi(1_A)} \xi_{\Phi(1_A)} (e_j, a_{ii}+b_{ii}) + \xi_{1_A} (\Phi(e_j), a_{ii}+b_{ii}) + \xi_{1_A}(e_j, \Phi(a_{ii})+\Phi(b_{ii})).
\end{align*}
Hence $\xi_{1_A}(e_j, t)=0$ which implies $t_{ij} = t_{ji} = t_{jj} = 0$. Now we get
$t = t_{ii}$.
For every $c_{ij} \in A_{ij}$, $i \neq j$, it follows from Claim \ref{claim6} that 
\begin{align*}
\sum_{\Phi(1_A)} & \xi_{\Phi(1_A)} (a_{ii}+b_{ii}, c_{ij}) + \xi_{1_A}(a_{ii}+b_{ii}, c_{ij}) + \xi_{1_A} (a_{ii}+b_{ii},\Phi(c_{ij}))\\
=& \ \Phi(\xi_{1_A} (a_{ii}+b_{ii}, c_{ij}))\\
=& \ \Phi(\xi_{1_A}(a_{ii},c_{ij})) + \Phi(\xi_{1_A}(b_{ii},c_{ij}))\\
=& \ \sum_{\Phi(1_A)} \xi_{\Phi(1_A)} (a_{ii}+b_{ii},c_{ij}) + \xi_{1_A} (\Phi(a_{ii})+\Phi(b_{ii}),c_{ij}) + \xi_{1_A} (a_{ii}+b_{ii}, \Phi(c_{ij})).
\end{align*}
Hence $\xi_{1_A}(t_{ii},c_{}ij)=0.$ Then $t_{ii}c_{ij} = 0$ for all $c_{ij} \in A_{ij}$, that is, $t_{ii}c \cdot e_j = 0$
for all $c \in A$. It follows from ($\spadesuit$) and ($\clubsuit$) that $t = t_{ii} = 0$, proving the claim. 

\begin{claim}\label{claim8}
$\Phi$ is additive.
\end{claim}
Let $a =\sum_{i,j=1}^{2} a_{ij}, \ b =\sum_{i,j=1}^{2} b_{ij} \in A$. By Claim \ref{claim5}, Claim \ref{claim6} and Claim \ref{claim7},
we have
\begin{align*}
\Phi(a + b) 
=& \Phi\left( \sum_{i,j=1}^{2} a_{ij} + \sum_{i,j=1}^{2} b_{ij} \right) = \Phi\left(\sum_{i,j=1}^{2}(a_{ij} + b_{ij})\right)\\
=& \ \sum_{i,j=1}^{2} \Phi(a_{ij} + b_{ij})\\
=& \ \sum_{i,j=1}^{2} \Phi(a_{ij}) +
\sum_{i,j=1}^{2} \Phi(b_{ij})\\
=& \ \Phi\left(\sum_{i,j=1}^{2} a_{ij}\right) +\Phi\left(\sum_{i,j=1}^{2} b_{ij}\right) = \Phi(a) + \Phi(b).
\end{align*}

\begin{claim}\label{claim9}
$\Phi(1_{A}) \bullet a = 0$ for all $a \in A$.
\end{claim}
Since $\Phi(1_A) \bullet 1_A = \Phi(1_A) + \Phi(1_A)^{*}$ is self-adjoint element, we have
\begin{align*}
2^{n-2} \Phi(1_A)
=& \ \Phi(\xi_{1_A}(1_A, 1_A))\\
=& \ \sum_{\Phi(1_A)} \xi_{\Phi(1_A)} (1_A,1_A) + \xi_{1_A}(\Phi(1_A), 1_A) + \xi_{1_A}(1_A,\Phi(1_A))\\
=& \ (n-1)2^{n-2}(\Phi(1_A)\bullet 1_A) + 2^{n-1}\Phi(1_A),
\end{align*}
which implies $\Phi(1_A) \bullet 1_A = 0$.
For every $a \in A$, it follows that
\begin{align*}
2^{n-1}\Phi(a)
=& \ \Phi(\xi_{1_A}(1_A,a))\\
=& \ \xi_{1_A}(\Phi(1_A), a) + \xi_{1_A}(1_A, \Phi(a))\\
=& \ 2^{n-2}\Phi(1_A)\bullet a + 2^{n-1} \Phi(a)
\end{align*}

Thus $\Phi(1_A) \bullet a = 0$.

\begin{claim}\label{claim10}
$\Phi(e_i) = \Phi(e_i)^{*}$, \ $i = 1, 2$.
\end{claim}

For every $a = a^{*} \in A$, since
\[
\xi_{1_A}(a,e_i) = \xi_{1_A}(e_i,a),
\]
we have
\[
\Phi(\xi_{1_A}(a,e_i)) = \Phi(\xi_{1_A}(e_i,a)).
\]
It follows from Claim \ref{claim9} that
\begin{align*}
\xi_{1_A}(\Phi(a),e_i) + \xi_{1_A}(a,\Phi(e_i)) = \xi_{1_A}(\Phi(e_i,a)) + \xi_{1_A}(e_i,\Phi(a))
\end{align*}
Then it is easy to compute that
$$a(\Phi(e_i) - \Phi(e_i)^{*}) = e_i(\Phi(a) - \Phi(a)^{*}).$$
We multiply above equation by $e_j$ from left side, it follows that
$$e_j \cdot a(\Phi(e_i) - \Phi(e_i)^{*}) = 0$$
for all $a = a^{*} \in A$. Thus for every $b \in A$, since $b = b_1 + ib_2$ with $b_1 = \frac{b+b^{*}}{2}$
and $b_2 = \frac{b-b^{*}}{2i}$ , we have $e_j \cdot b(\Phi(e_i) - \Phi(e_i)^{*}) = 0$
for all $b \in A$. Thus
$$(\Phi(e_i)^{*} - \Phi(e_i))b \cdot e_j = 0$$
for all $b \in A$. It follows from ($\spadesuit$) and ($\clubsuit$) that $\Phi(e_i) = \Phi(e_i)^{*}$.

\begin{claim}\label{claim11}
$\Phi(e_i) = e_1\Phi(e_i)e_2 + e_2\Phi(e_i)e_1$ \ $(i = 1, 2)$ and $\Phi(1_A) = 0.$
\end{claim}

On the one hand, by Claim \ref{claim9} and Claim \ref{claim10}, we have
\[
2^{n-1}\Phi(e_i) = \Phi(\xi_{1_A}(e_i,e_i)) = \xi_{1_A}(\Phi(e_i),e_i) + \xi_{1_A}(e_i, \Phi(e_i)) = 2^{n-1} (\Phi(e_i)e_i + e_i\Phi(e_i)).
\]
Multiplying both sides of the above equation by $e_j$ $\left(1 \leq i \neq j \leq 2 \right)$ from the left
and right respectively, we obtain that $e_j \Phi(e_i) e_j = 0$. Similarly, multiplying both
sides of the above equation by $e_i$ from the left and right respectively, we have
that $e_i\Phi(e_i)e_i = 0$. Hence $\Phi(e_i) = e_1\Phi(e_i)e_2 + e_2\Phi(e_i)e_1$, $i = 1, 2$.
On the other hand, we have
\begin{align*}
0 
=& \ \Phi(\xi_{1_A}(e_1,e_2)) = \xi_{1_A}(\Phi(e_1),e_2) + \xi_{1_A}(e_1,\Phi(e_2))\\
=& \  2^{n-2} (e_1\Phi(e_1)e_2 + e_2\Phi(e_1)e_1 + e_1\Phi(e_2)e_2 + e_2\Phi(e_2)e_1).
\end{align*}

Now, by the additivity of $\Phi$, we can get that
$$\Phi(1_A) = \Phi(e_1)+\Phi(e_2) = e_1\Phi(e_1)e_2 +e_2\Phi(e_1)e_1 +e_1\Phi(e_2)e_2 +e_2\Phi(e_2)e_1 = 0.$$

\begin{claim}\label{claim12}
For all $a, b \in A$, we have $\Phi(a \bullet b) = \Phi(a) \bullet b + a \bullet \Phi(b)$.
\end{claim}
For every $a, b \in A$, it follows from $\Phi(1_A) = 0$ that
\[
2^{n-2}\Phi(a \bullet b)= \Phi(\xi_{1_A}(a,b)) = \xi_{1_A}(\Phi(a),b) + \xi_{1_A}(a,\Phi(b)) = 2^{n-2}(\Phi(a)\bullet b + a \bullet \Phi(b)).
\]

Hence $\Phi(a \bullet b) = \Phi(a) \bullet b + a \bullet \Phi(b)$.

\begin{claim}\label{claim13}
For all $a \in A$, $\Phi(a^{*}) = \Phi(a)^{*}.$
\end{claim}

For every $a \in A$, by Claims \ref{claim8}, \ref{claim11}, \ref{claim12}, we have
$\Phi(a) + \Phi(a^{*}) = \Phi(a \bullet 1_A) = \Phi(a) \bullet 1_A = \Phi(a) + \Phi(a)^{*}$.
Hence $\Phi(a^{*}) = \Phi(a)^{*}$.

Now let us define a map $\phi : A \rightarrow A$ by $\phi(a) = \Phi(a) - \Xi_{y,z}(a)$ for all $a \in A$, where $\Xi_{y,z}$ is $*$-inner derivation on $A$, with $y = e_1 \Phi(e_1) e_2 + e_2 \Phi(e_1)e_1$ and $z = e_1$. It is easy to see the following properties.

\begin{claim}\label{claim14}

\begin{enumerate}
	\item [\it (1)] For all $a,b \in A$, $\phi(a \bullet b) = \phi(a) \bullet b + a \bullet \phi(b)$;
	\item [\it (2)] $\phi(e_i) = 0$, $i = 1, 2$;
	\item [\it (3)] $\phi$ is additive;
	\item [\it (4)] For all $a \in A$, $\phi(a^{*}) = \phi(a)^{*}$;
	\item [\it (5)] $\phi$ is an additive derivation if and only if $\Phi$ is an additive derivation.
\end{enumerate}

\end{claim}

\begin{claim}\label{claim15}
$\phi(a_{ij}) \in A_{ij}$, $i, j = 1, 2$.
\end{claim}

Let $a_{ij} \in A_{ij}$, $1 \leq i \neq j \leq 2$. On the one hand, it follows from $\phi(e_i) = 0$ that
$$\phi(a_{ij}) = \phi(e_i \bullet a_{ij}) = e_i \bullet \phi(a_{ij}) = e_i\phi(a_{ij}) + \phi(a_{ij})e_i.$$
Hence $e_i\phi(a_{ij})e_i = e_j\phi(a_{ij})e_j = 0$. On the other hand, we have
$$0 = \phi(a_{ij} \bullet e_i) = \phi(a_{ij}) \bullet e_i = \phi(a_{ij})e_i + e_i\Phi(a_{ij})^{*}.$$
Hence $e_j\phi(a_{ij})e_i = 0$. Thus $\phi(a_{ij}) = e_i\phi(a_{ij})e_j \in A_{ij}$, $1 \leq i \neq j \leq 2$.
Let $a_{ii} \in A_{ii}$, $i = 1, 2$. Then
$$0 = \phi(e_j \bullet a_{ii}) = e_j \bullet \phi(a_{ii}) = e_j\phi(a_{ii}) + \phi(a_{ii})e_j.$$
Hence $e_i\Phi(a_{ii})e_j = e_j\Phi(a_{ii})e_i = e_j\Phi(a_{ii})e_j = 0$. Now we get $\Phi(a_{ii}) =
e_i\Phi(a_{ii})e_i \in A_{ii}$, $i = 1, 2$.

\begin{claim}\label{claim16}
Let $a_{ii},b_{ii} \in A_{ii}$ and $a_{ij}, b_{ij} \in A_{ij}$, $1 \leq i \neq j \leq 2$. Then
$$\phi(a_{ii}b_{ii}) = \phi(a_{ii})b_{ii} + a_{ii}\phi(b_{ii}), \  \phi(a_{ii}b_{ij}) = \phi(a_{ii})b_{ij} + a_{ii}\phi(b_{ij}),$$
$$\phi(a_{ij}b_{ji}) = \phi(a_{ij})b_{ji} + a_{ij}\phi(b_{ji}), \  \phi(a_{ij}b_{jj}) = \phi(a_{ij})b_{jj} + a_{ij}\phi(b_{jj}),$$
$$\phi(a_{ij}b_{ij}) = \phi(a_{ij})b_{ij} + a_{ij}\phi(b_{ij}).$$

\end{claim}

It follows from Claim \ref{claim15} that
\[
\phi(a_{ii}b_{ij}) = \phi(a_{ii} \bullet b_{ij}) = \phi(a_{ii}) \bullet b_{ij} + a_{ii} \bullet \phi(b_{ij}) = \phi(a_{ii})b_{ij} + a_{ii}\phi(b_{ij}),
\]

that is $\phi(a_{ii}b_{ij}) = \phi(a_{ii})b_{ij} + a_{ii}\phi(b_{ij})$. By flexibility of alternative algebras and
for any $c_{ij} \in A_{ij}$, we have
\begin{align*}
\phi(a_{ii}b_{ii})\cdot c_{ij} + a_{ii}b_{ii} \cdot \phi(c_{ij}) 
=& \ \phi(a_{ii}b_{ii} \cdot c_{ij}) = \phi(a_{ii}) \cdot b_{ii}c_{ij} + a_{ii} \cdot \phi(b_{ii}c_{ij}) \\
=& \ \phi(a_{ii})\cdot b_{ii}c_{ij} + a_{ii} \cdot \phi(b_{ii})c_{ij} + a_{ii} \cdot b_{ii}\phi(c_{ij})\\
=& \ \phi(a_{ii}) b_{ii}\cdot c_{ij} + a_{ii}  \phi(b_{ii})\cdot c_{ij} + a_{ii}  b_{ii}\cdot \phi(c_{ij}).
\end{align*}

Then $(\phi(a_{ii}b_{ii})-\phi(a_{ii})b_{ii}-a_{ii}\phi(b_{ii}))c_{ij} = 0$ for any $c_{ij} \in A_{ij}$. It follows from
($\spadesuit$) and ($\clubsuit$) that $\phi(a_{ii}b_{ii}) = \phi(a_{ii})b_{ii} + a_{ii}\phi(b_{ii}).$
It follows from Claim \ref{claim15} that
\begin{align*}
\phi(a_{ij}b_{ji}) + \phi(b_{ji}a_{ij}^{*}) = \phi(a_{ij} \bullet b_{ji})
=& \ \phi(a_{ij}) \bullet b_{ji} + a_{ij} \bullet \phi(b_{ji})\\
=& \ \phi(a_{ij})b_{ji} + a_{ij}\phi(b_{ji}) + b_{ji}\phi(a_{ij})^{*} + \phi(b_{ji})a_{ij}^{*},
\end{align*}
hence $\phi(a_{ij}b_{ji}) = \phi(a_{ij})b_{ji} + a_{ij}\phi(b_{ji})$.
Again by flexibility of alternative algebras and for any $c_{ji} \in A_{ji}$ \ ($i \neq j$), we have
\begin{align*}
c_{ji}\cdot \phi(a_{ij}b_{jj}) + \phi(c_{ji})\cdot a_{ij}b_{jj} 
=& \ \phi(c_{ji}\cdot a_{ij}b_{jj}) = \phi(c_{ji} a_{ij}\cdot b_{jj}) = \phi(c_{ji}a_{ij})b_{jj} + c_{ji}a_{ij}\phi(b_{jj})\\
=& \ \phi(c_{ji})a_{ij} \cdot b_{jj} + c_{ji}\phi(a_{ij}) \cdot b_{jj} + c_{ji}a_{ij} \cdot \phi(b_{jj})\\ 
=& \ \phi(c_{ji})\cdot a_{ij}  b_{jj} + c_{ji}\cdot \phi(a_{ij})  b_{jj} + c_{ji}\cdot a_{ij}  \phi(b_{jj}).
\end{align*}
Hence $c_{ji}(\phi(a_{ij}b_{jj}) - \phi(a_{ij})b_{jj} - a_{ij}\phi(b_{jj})) = 0$ for any $c_{ji} \in A_{ji}$. It follows
from ($\spadesuit$) and ($\clubsuit$) that $\phi(a_{ij}b_{jj}) = \phi(a_{ij})b_{jj} + a_{ij}\phi(b_{jj})$.
Since $\phi$ is a $*$-Jordan derivation we have
\begin{eqnarray*}\label{ident}
&&\phi(a \bullet (bc + cb^{*})) = \phi(a) \bullet (bc + cb^{*}) + a \bullet \phi(bc + cb^{*})
\\&=& \phi(a) (bc + cb^{*}) + (bc + cb^{*}) \phi(a)^{*} + a \phi(bc + cb^{*}) + \phi(bc + cb^{*}) a^{*}.
\end{eqnarray*}
Tanking $a= e_j$, $b = a_{ij}$ and $c = b_{ij}$ in (\ref{ident}) we get
\begin{align*}
\phi(a_{ij}b_{ij}) 
=& \ e_j \phi(a_{ij}b_{ij} + b_{ij}a_{ij}^{*}) + \phi(a_{ij}b_{ij} + b_{ij}a_{ij}^{*}) e_j \\
=& \ e_j \cdot \phi(a_{ij})b_{ij} + e_j \cdot a_{ij}\phi(b_{ij}) +  e_j \cdot \phi(b_{ij})a_{ij}^{*} + e_j \cdot b_{ij}\phi(a_{ij}^{*})
\\
=& \ \phi(a_{ij})b_{ij} + a_{ij}\phi(b_{ij}). 
\end{align*}

\begin{claim}\label{claim17}
$\Phi(ab) = \phi(a)b + a\phi(b)$ for all $a,b \in A.$
\end{claim}
 
Write $a = \sum_{i,j = 1}^{2} a_{ij}, b = \sum_{i,j = 1}^{2} b_{ij} \in A$. Then $ab = a_{11}b_{11} + a_{11}b_{12} +
a_{12}b_{12} + a_{12}b_{21} + a_{12}b_{22} + a_{21}b_{11} + a_{21}b_{12} + a_{21}b_{21} + a_{22}b_{21} + a_{22}b_{22}$. It follows from Claim \ref{claim16} and the additivity of $\phi$ that $\phi(ab) = \phi(a)b +a\phi(b)$. By the definition of $\phi$, we obtain that $\Phi(ab) = \Phi(a)b + a\Phi(b)$.
Now, by Claim \ref{claim8}, Claim \ref{claim13} and Claim \ref{claim16}, we have proved that $\Phi$ is an additive $*$-derivation. This completes the proof.

\section{Application}

A complete normed alternative complex $*$-algebra $A$ is called of alternative
$C^{*}$-algebra if it satisfies the condition: $\left\|a^{*}a\right\| = \left\|a\right\|^2$, for all elements $a \in A$. 
A non-zero element $p \in A$ is called a projection if it is self-adjoint and verifies the
condition $p^2 = p$. Alternative $C^{*}$-algebras are non-associative generalizations of
$C^{*}$-algebras and appear in various areas in Mathematics (see more details in the references \cite{Miguel1} and \cite{Miguel2}). 
An alternative $C^{*}$-algebra $A$ is called of alternative $W^{*}$-algebra if it is a dual Banach space and a prime alternative $W^{*}$-algebra is called
alternative $W^{*}$-factor. It is well known that non-zero alternative $W^{*}$-algebras are unital and that an alternative $W^{*}$-algebra $A$ is a factor if and only if its center is equal to $\mathbb{C}1_A$, where $1_A$ is the unit of $A$.

\begin{theorem}
Let $A$ be an alternative $W^{*}$-factor. 
If a map $\Phi : A \rightarrow A$ satisfies
$$\Phi(a_1 \bullet a_2 \bullet \cdots \bullet a_n) =
\sum_{k=1}^n a_1 \bullet \cdots \bullet a_{k-1} \bullet \Phi(a_k) \bullet a_{k+1} \bullet \cdots \bullet a_n$$
for all $a_{n-1}, a_n \in A$ and $a_i = 1_A$ for all $i \in \left\{1, 2, \ldots, n-2\right\}$, then $\Phi$ is an additive $*$-derivation.  
\end{theorem}

\begin{corollary}
Let $A$ be an alternative $W^{*}$-factor. 
Then $\Phi$ is a nonlinear $*$-Jordan-type derivation on $A$ if and only if $\Phi$ is an additive
$*$-derivation.
  
\end{corollary}

\section{Conflict of interest}
On behalf of all authors, the corresponding author states that there is no conflict of interest.


\begin{thebibliography}{99}
  
\bibitem{BreFo}
Bre$\check{s}$ar M., Fo$\check{s}$ner
On ring with involution equipped with some new product, 
Publ. Math. Debrecen, 57 (2000) 121--134.  
  
  
\bibitem{mff16}
    Da Motta Ferreira, J. C., Ferreira, B. L. M., Additivity of -Multiplicative Maps on Alternative Rings,  Communications In Algebra (Online), 44 (2016), 1557-1568, 2016.	
	
	
\bibitem{Dai} 
Dai L., Lu F., 
    Nonlinear maps preserving Jordan $*$-products, 
    J. Math. Anal. Appl.,  409 (2014), 180--188.
    
\bibitem{Darvish}
	Darvish V. ,  Nouri M.,  Razeghi M., Taghavi A.
	Nonlinear $*$-Jordan triple
derivations on prime $*$-algebras,
Rocky Mountain J. Math., 50 (2020) 543--549.
  
\bibitem{f19}
    Ferreira, B. L. M.,  Additivity of elementary maps on alternative rings, Algebra and Discrete Mathematics, 28 (2019),  94--106.  

\bibitem{FerreiraCosta}
Ferreira B. L. M., Costa B. T., $*$-Jordan-type maps on $C^{*}$-algebras,
	arXiv:2005.11430 [math.OA] (2020)  
	
	
\bibitem{Fer} 
    Ferreira J. C. M., Ferreira B. L. M., 
    Additivity of $n$-multiplicative maps on alternative rings, 
    Communications in Algebra, 44 (2016), 1557--1568.	

\bibitem{fgf20}
    Ferreira, B. L. M., Guzzo, H., Ferreira, R. N.,  An Approach Between the Multiplicative and Additive Structure of a Jordan Ring, Bulletin of the Iranian Mathematical Society, 1 (2020), 1.

\bibitem{posd} 
  Ferreira B. L. M., Guzzo H., Lie Maps on alternative rings, Bollettino dell'Unione Matematica Italiana, 13 (2020), 2, 181--192.
    
\bibitem{fgk20} 
    Ferreira B. L. M., Guzzo H., Kaygorodov I., Lie maps on alternative rings preserving idempotents, Colloquium Mathematicum, (2021), DOI: 10.4064/cm8195-10-2020.
	
\bibitem{Miguel1} 
    Garc\'ia M. C., Palacios \'A. R., 
    Non-associative normed algebras. Vol. 1. The Vidav-Palmer and Gelfand-Naimark theorems. Encyclopedia of Mathematics and its Applications, 154. Cambridge University Press, Cambridge, 2014. xxii+712 pp.


\bibitem{Miguel2}  Garc\'ia M. C., Palacios \'A. R., 
    Non-associative normed algebras. Vol. 2. Representation theory and the Zel'manov approach. Encyclopedia of Mathematics and its Applications, 167. Cambridge University Press, Cambridge, 2018. xxvii+729 pp.

\bibitem{He} 
    Hentzel I.,  Kleinfeld E.,  Smith H., 
    Alternative rings with idempotent, 
    Journal of Algebra, 64 (1980), 325--335.
	
\bibitem{Huo} 
Huo D., Zheng B. and Liu H., Nonlinear maps preserving Jordan triple $\eta$-$*$-products, J. Math. Anal. Appl.,  430 (2015), 830--844.
		
\bibitem{LiLuFang2}	
 Li C., Lu F., Fang X., Nonlinear $\epsilon$-Jordan $*$-derivations on von Neumann algebras,  Linear Multilinear Algebra., 62 (2014) 466--473.
		
\bibitem{LiLuFa} 
Li C., Lu F. and Fang X.,
    Mappings preserving new product $XY +Y X^{*}$ on
factor von Neumann algebras, 
    Linear Algebra Appl.,  438 (2013), 2339--2345.
    
		
\bibitem{LiLu}	
	Li C., Lu F.,	
	 Nonlinear maps preserving the Jordan triple $*$-product on von
Neumann algebras,
 Ann. Funct. Anal., 7 (2016), 496--507.

\bibitem{LiLu2}
Li C., Lu F.,
Nonlinear Maps Preserving the Jordan Triple $1$-$*$-product on
von Neumann algebras, 
Complex Anal. Oper. Theory., 11 (2017), 109--117.

\bibitem{Molnar}
Moln$\acute{a}$r L.,
A condition for a subspace of $B(H)$ to be an ideal,
Linear Algebra Appl., 235 (1996) 229--234.
		


\bibitem{Sem1}
$\check{S}$emrl P.,
Quadratic functionals and Jordan $*$-derivations,
Studia Math., 97 (1991) 157--165.


\bibitem{Sem3}
$\check{S}$emrl P.,
On Jordan $*$-derivations and an application,
Colloq. Math., 59 (1990) 241--251.









\bibitem{TaghaRo}
Taghavi A., Rohi H. and Darvish V., 
Non-linear $*$-Jordan derivations on von
Neumann algebras, 
Linear Multilinear Algebra, 64 (2016), 426--439.
 
\bibitem{Zhang}
Zhang F., 
Nonlinear skew Jordan derivable maps on factor von Neumann
algebras,
Linear Multilinear Algebra, 64 (2016) 2090--2103.

\bibitem{ZhaoLi}		
	 Zhao	F., Li C., 
Nonlinear $*$-Jordan triple derivations on von Neumann algebras,
Math. Slovaca, 68 (2018) 163--170.

\bibitem{ZhaLi}	
Zhao F., Li C.,		
 Nonlinear maps preserving the Jordan triple $*$-product between factors,
 Indag. Math., 29 (2018) 619--627.


 

    
    
    


    


    

    

    







    
    
    
 
 
     



\end{thebibliography}
\end{document}